\documentclass[12pt]{article}
\usepackage{amsmath,amsthm,amssymb,amscd}
\usepackage{latexsym}
\newtheorem{thm}{Theorem}[section]

 \newtheorem{prop}[thm]{Proposition}
 \theoremstyle{definition}
 
 \theoremstyle{remark}

  \numberwithin{equation}{section}

\makeatletter
\newcommand*\owedge{\mathpalette\@owedge\relax}
\newcommand*\@owedge[1]{%
  \mathbin{%
    \ooalign{%
      $#1\m@th\bigcirc$\cr
      \hidewidth$#1\m@th\wedge$\hidewidth\cr
    }%
  }%
}
\makeatother

\begin{document}

\title
{Open manifolds  with uniformly positive isotropic curvature}

\author{ Hong Huang}
\date{}
\maketitle

\begin{abstract}
 We prove the following result: Let $(M,g_0)$ be a complete noncompact
manifold of dimension $n\geq 12$ with isotropic curvature bounded below by a positive constant, with scalar curvature  bounded above, and with injectivity radius bounded below. Then  there is a finite collection
$\mathcal{X}$ of spherical $n$-manifolds and manifolds of the form $\mathbb{S}^{n-1} \times
\mathbb{R} /G$, where $G$ is a   discrete subgroup
of the isometry group of the round cylinder $\mathbb{S}^{n-1}\times
\mathbb{R}$, such that $M$ is diffeomorphic to a (possible infinite)
connected sum of  members of $\mathcal{X}$.
This extends a recent work of Huang.  The proof uses Ricci flow with surgery on open orbifolds with isolated singularities.


\end{abstract}


\section {Introduction}

In a previous paper \cite{Hu23} we get a topological classification of compact manifolds of dimension $n\geq 12$ with positive isotropic curvature, using a work of Brendle \cite{B19} and a previous work of Huang \cite{Hu22}. In this note we extend the main result in \cite{Hu23} to the case of open (i.e. noncompact and complete) manifolds  of dimension $n\geq 12$ with uniformly positive isotropic curvature and bounded geometry (in the sense of having bounded sectional curvature and injectivity radius lower bound).

Recall that in our previous work \cite{Hu15}, we say that a manifold of dimension $n\geq 4$ has  uniformly positive isotropic curvature if its isotropic curvature is bounded below by a positive constant. (Note that in  \cite{B19}, \cite{BN1}, and \cite{BN2}, a manifold of dimension $n\geq 5$ is said to have uniformly positive isotropic curvature if $\text{Rm}-\theta R \hspace*{1mm} \text{id} \owedge \text{id} \in \text{PIC}$ for some positive constant $\theta$, where $\text{PIC}$ denotes the cone of algebraic curvature tensors in dimension $n$ with nonnegative isotropic curvature, and there is also a definition of  uniformly positive isotropic curvature in dimension $n=4$ in \cite{BN1} and \cite{BN2}.  Their definition is slightly different from that of ours.)

Let's also recall the notion of (possibly infinite) connected sum,
compare for example p. 204-205 of \cite{M07},  \cite{BBM} and \cite{Hu17}. Let  $\mathcal{X}$ be a set of $n$-manifolds. A
$n$-manifold $M$ is said to be a connected sum of members of
$\mathcal{X}$ if there exist a (not necessarily locally finite)  graph  and a map
$v\mapsto X_v$ which associates a member in $\mathcal{X}$ to each vertex of this graph, such that if for each edge of the graph, supposing  it connects the vertices $v_1$ and $v_2$ (here we allow $v_1=v_2$), we do a connected sum  of $X_{v_1}$ and $X_{v_2}$ (or a self connected sum of $X_{v_1}$ when $v_1=v_2$), we obtain a $n$-manifold
diffeomorphic to $M$.  Here we also allow more than one edge connecting two vertices (or one vertex to itself).  We can also extend this notion to the case of orbifolds.  For the definition of  connected sum of two smooth manifolds/orbifolds or a self connected sum, see for example,  \cite{B}, \cite{BJ},  \cite{K}, \cite{M07}, \cite{Mu}, and \cite{Hu22}.

Our main result is

    \begin{thm} \label{thm 1.1} Let $(M,g_0)$ be a complete noncompact
manifold of dimension $n\geq 12$ with isotropic curvature bounded below by a positive constant, with scalar curvature  bounded above, and with injectivity radius bounded below. Then  there is a finite collection
$\mathcal{X}$ of spherical $n$-manifolds and manifolds of the form $\mathbb{S}^{n-1} \times
\mathbb{R} /G$, where $G$ is a   discrete subgroup
of the isometry group of the round cylinder $\mathbb{S}^{n-1}\times
\mathbb{R}$, such that $M$ is diffeomorphic to a (possible infinite)
connected sum of  members of $\mathcal{X}$.
\end{thm}

 \noindent (By  \cite{MW} the converse of Theorem 1.1 is also true.)  Compare Theorem 1.1 in \cite{Hu15} and Theorem 1.1 in \cite{Hu23}.

In the process of proof of Theorem 1.1 we get a slightly more general result.

\begin{thm} \label{thm 1.2} \ \  Let $(\mathcal{O},g_0)$ be a complete noncompact
orbifold of dimension $n\geq 12$ with at most isolated singularities, with   isotropic curvature bounded below by a positive constant, with scalar curvature  bounded above, and with $\text{Vol} \hspace*{1mm} B_1(x)\geq v>0$ for every $x \in |\mathcal{O}|$. Then   there is a finite collection
$\mathcal{X}$ of spherical $n$-orbifolds  with  at most isolated singularities  such that  $\mathcal{O}$ is diffeomorphic to  an infinite
connected sum of  members of $\mathcal{X}$.
\end{thm}

\noindent (By an obvious extension of Theorem 1.1 in \cite{MW} to the case of orbifold connected sum (cf. also \cite{CH}), certain version of the converse of Theorem 1.2 is also true.) Compare Theorem 1.3 in \cite{Hu15}  and Theorem 1.3 in \cite{Hu23}.

Note that Brendle's curvature pinching estimates for Ricci flow on compact manifolds of dimension $n\geq 12$ with positive isotropic curvature (see Theorem 1.2 and Corollary 1.3 in \cite{B19}) extend to the case of complete (noncompact)
manifolds/orbifolds of dimension $n\geq 12$ with isotropic curvature bounded below by a positive constant, and with scalar curvature  bounded above. With this and  \cite{Hu22} and \cite{Hu23}  at hand, the main work in the present note is to argue the finiteness of the collection $\mathcal{X}$ appearing in Theorems \ref{thm 1.1} and \ref{thm 1.2}.  By the way, our previous work \cite{Hu15} is somewhat sloppy in some places.

 In Section 2 we  describe the canonical neighborhood structure of orbifold ancient
 $\kappa$-solutions. In Section 3 we  construct  an $(r, \delta)$-surgical solution to the Ricci flow starting with a complete noncompact
orbifold of dimension $n\geq 12$ with at most isolated singularities,  with   isotropic curvature bounded below by a positive constant, with scalar curvature  bounded above, and with $\text{Vol} \hspace*{1mm} B_1(x)\geq v>0$ for every $x \in |\mathcal{O}|$. Finally, in Section 4, we prove
Theorems 1.1 and  1.2.  We will follow the conventions and notation in \cite{Hu23} in general.

\section{Orbifold ancient $\kappa$-solutions}

The following result is an extension of Corollary 6.22 in \cite{B19}, and a refinement of Proposition 2.7 in \cite{Hu23}.
\begin{prop}\label{prop 2.1} (cf. \cite{P2}, Proposition 6.21 and Corollaries 6.20 and 6.22 in \cite{B19}, and Proposition 2.7 in \cite{Hu23}) \ \
Given   a small positive constant $\varepsilon$  and a constant $\theta >0$, there exist positive constants $C_1=C_1(n, \theta,\varepsilon)$, $C_2=C_2(n,\theta,\varepsilon)$ and $C_3=C_3(n,\theta)$ with the following property:
Suppose that $(\mathcal{O},g(t))$ is an orbifold ancient $\kappa$-solution of dimension $n\geq 5$ with at most isolated singularities which satisfies
$Rm-\theta R \hspace*{1mm} \text{id} \owedge \text{id} \in PIC$. Then  for each space-time point $(x_0,t_0)$, there
is  an open subset $U$ of $\mathcal{O}$ with
$B(x_0,t_0,C_1^{-1}R(x_0,t_0)^{-\frac{1}{2}})\subset U \subset B(x_0,t_0, C_1R(x_0,t_0)^{-\frac{1}{2}})$ and with  $C_2^{-1}R(x_0,t_0) \leq R(x,t_0) \leq C_2R(x_0,t_0)$ for any $x\in U$, which falls into
one of the following three categories:

(a) $U$ is a strong $\varepsilon$-neck  centered at $(x_0,t_0)$,

(b) $U$ is an $\varepsilon$-cap centered at $(x_0,t_0)$,  and

(c)  $U$ is  compact and strictly PIC2;

\noindent moreover, the scalar curvature in $U$ at time $t_0$ satisfies
the derivative estimates
\begin{equation*}
|\nabla R|\leq C_3 R^{\frac{3}{2}} \hspace*{8mm} \mbox{and} \hspace*{8mm}
|\frac{\partial R}{\partial t}|\leq C_3 R^2,
\end{equation*}
and in case (c) the sectional curvature in $U$ at time $t_0$ is greater than $C_2^{-1}R(x_0,t_0)$.
\end{prop}

\noindent {\bf Proof}. \ \ Compare the proof of Proposition 6.21 and Corollary 6.22 in \cite{B19}. By Proposition 2.7  in \cite{Hu23},  we only need to consider the case that $\mathcal{O}$ is compact. Fix $t_0=0$. Suppose that the result is not true in this case. Then there is a sequence of compact ancient orbifold $\kappa_j$-solutions $(\mathcal{O}^{(j)},g^{(j)}(t))$ of dimension $n\geq 5$ with $\text{Rm}-\theta R \hspace*{1mm} \text{id} \owedge \text{id} \in \text{PIC}$ and a sequence of points $x_j\in \mathcal{O}^{(j)}$ with the following property: For each $j$, there does not exist a neighborhood $U$ of $x_j$ with $B_{g^{(j)}(0)}(x_j, j^{-1}R(x_j,0)^{-\frac{1}{2}}) \subset U \subset B_{g^{(j)}(0)}(x_j, jR(x_j,0)^{-\frac{1}{2}})$ and $j^{-1}R(x_j,0)\leq R(x,0) \leq jR(x_j,0)$ for all $x\in U$, and such that $U$  is either a strong  $\varepsilon$-neck  centered at $(x_j,0)$, or an $\varepsilon$-cap centered at $(x_j,0)$, or is  compact,  strictly PIC2, and at time $0$ has  sectional curvature   greater than $j^{-1}R(x_j, 0)$.  By rescaling we can assume that $R(x_j,0)=1$ for all $j$.

By using an orbifold version of Corollary 6.7 in \cite{B19} we see that each $(\mathcal{O}^{(j)},g^{(j)}(t))$ is strictly PIC2. So $\mathcal{O}^{(j)}$ is diffeomorphic to a spherical orbifold $\mathbb{S}^n//G_j$ with at most isolated singularities, where $G_j$ is a  finite subgroup of $O(n+1)$ acting freely on $\mathbb{S}^n$, or a finite subgroup of $O(n)\times \mathbb{Z}_2$  by  Lemma 5.2 in \cite{CTZ}, Lemma  3.2 in \cite{Hu22}, and Corollary 2.4 in Chapter VI of \cite{B72}. We pull back  $(\mathcal{O}^{(j)},g^{(j)}(t))$ to a $G_j$-invariant ancient $\kappa_j$-solution  $(M^{(j)},\tilde{g}^{(j)}(t))$ via the universal covering maps $\pi_j:M^{(j)}\rightarrow \mathcal{O}^{(j)}$, where each $M^{(j)}$ is diffeomorphic to $\mathbb{S}^n$. Choose $\tilde{x}_j \in \pi_j^{-1}(x_j)$. By  Theorem 6.19 in \cite{B19} $(M^{(j)},\tilde{g}^{(j)}(t))$ is $\kappa_0$-noncollapsed for some positive constant
$\kappa_0$ independent of $j$.  By \cite{BN1}  each $(M^{(j)},\tilde{g}^{(j)}(t))$ is  $O(n)\times \mathbb{Z}_2$-invariant.

By Corollary 6.15 in \cite{B19}, and Lemma 3.6 and Proposition 3.7 in \cite{BHZ},  after passing to a subsequence (we will use the same notation for the  sequence itself and its subsequence), the sequence $(M^{(j)},\tilde{g}^{(j)}(t), (\tilde{x}_j,0))$ smoothly $O(n)$-equivariantly converges to an  $O(n)$-invariant ancient $\kappa_0$-solution $(M^{(\infty)},\tilde{g}^{(\infty)}(t), (\tilde{x}_\infty,0))$ with  $\text{Rm}-\theta R \hspace*{1mm} \text{id} \owedge \text{id} \in \text{PIC}$.

If $M^{(\infty)}$ is compact, then the diameter of $(M^{(j)},\tilde{g}^{(j)}(0))$ has a uniform upper bound independent of $j$. So for $j$ large enough $\tilde{U}_j:=M^{(j)}$ is a neighborhood of $\tilde{x}_j$ with $B_{\tilde{g}^{(j)}(0)}(\tilde{x}_j, j^{-1}) \subset \tilde{U}_j \subset B_{\tilde{g}^{(j)}(0)}(\tilde{x}_j, j)$ and $j^{-1}\leq R(\tilde{x},0) \leq j$ for all $\tilde{x}\in \tilde{U}_j$; moreover  at time $0$, $\tilde{U}_j$ has  sectional curvature  greater than $j^{-1}$.  Using the fact that $\pi_j|_{B(\tilde{x}_j,0,r)}:B(\tilde{x}_j,0,r)\rightarrow  B(x_j,0,r)$ is surjective for any $r>0$,  we get that for $j$ large enough,
$B_{g^{(j)}(0)}(x_j, j^{-1}) \subset \mathcal{O}^{(j)} \subset B_{g^{(j)}(0)}(x_j, j)$,  and $j^{-1}\leq R(x,0) \leq j$ for all $x\in \mathcal{O}^{(j)}$; moreover, $\mathcal{O}^{(j)}$ is  compact,  strictly PIC2,  and at time $0$ has sectional curvature    greater than $j^{-1}$. A contradiction.

If $M^{(\infty)}$ is noncompact,  by Proposition 2.2 in \cite{BHZ} $M^{(\infty)}$  is diffeomorphic to $\mathbb{S}^{n-1} \times \mathbb{R}$ or $\mathbb{R}^n$. In this case, when $j$ is sufficiently large, $(M^{(j)},\tilde{g}^{(j)}(t))$ must coincide with Perelman's solution (which is originally constructed in Section 1.4 of \cite{P2} for $n=3$, but the construction can be extended to the higher dimensional case) up to diffeomorphisms,
translations in time, and parabolic rescalings by Theorem 1.4 in \cite{BN1}. Note that the isometry group of Perelman's solution (in dimension $n$) is exactly  $O(n)\times \mathbb{Z}_2$. So in this case, when $j$ is sufficiently large, $G_j$ must be a finite subgroup of $O(n)\times \mathbb{Z}_2$. If $M^{(\infty)}$  is diffeomorphic to $\mathbb{S}^{n-1} \times \mathbb{R}$,  $(M^{(\infty)},\tilde{g}^{(\infty)}(t))$ can not be strictly $\text{PIC2}$, then by Corollary 6.7 in \cite{B19} we see that $(M^{(\infty)},\tilde{g}^{(\infty)}(t))$ is a shrinking round cylinder. So for any $\eta>0$, when $j$ is sufficiently large, there exists an $O(n)$-invariant strong $\eta$-neck $ \tilde{U}_j \subset M^{(j)}$ centered at $(\tilde{x}_j, 0)$ with $B_{\tilde{g}^{(j)}(0)}(\tilde{x}_j, \frac{1}{2}\eta) \subset \tilde{U}_j \subset B_{\tilde{g}^{(j)}(0)}(\tilde{x}_j, 2\eta^{-1})$ and $\frac{1}{2}\leq R(\tilde{x},0) \leq 2$ for all $\tilde{x}\in \tilde{U}_j$.  Then for suitably chosen $\eta$ and large $j$, certain subset $U$ of $\pi_j(\tilde{U}_j)$ is (a time $0$-slice of) a strong $\varepsilon$-neck or an $\varepsilon$-cap centered at $(x_j,0)$,  (for the cap case we  need to use  the argument in the proof of Proposition  2.4 in \cite{Hu23} to attain the extra property required in the definition of an $\varepsilon$-cap in \cite{Hu23},)  with $B_{g^{(j)}(0)}(x_j, \frac{1}{2}\eta) \subset U \subset B_{g^{(j)}(0)}(x_j, 2\eta^{-1})$, and $\frac{1}{2} \leq R(x,0) \leq 2$ for all $x\in U$.  Thus we get a contradiction.   If $M^{(\infty)}$  is diffeomorphic to $\mathbb{R}^n$, by the proof of Theorem 6.18 in \cite{B19} and Proposition 2.6 in \cite{Hu23} and the fact that the sequence $(M^{(j)},\tilde{g}^{(j)}(t), (\tilde{x}_j,0))$ smoothly $O(n)$-equivariantly converges to  $(M^{(\infty)},\tilde{g}^{(\infty)}(t), (\tilde{x}_\infty,0))$,   when $j$ is sufficiently large, there exists an $O(n)$-invariant strong $\varepsilon$-neck $ \tilde{U}_j$ or an $O(n)$-invariant $\varepsilon$-cap $ \tilde{U}_j$ centered at $(\tilde{x}_j, 0)$ with $B_{\tilde{g}^{(j)}(0)}(\tilde{x}_j, (2C_1)^{-1}) \subset \tilde{U}_j \subset B_{\tilde{g}^{(j)}(0)}(\tilde{x}_j, 2C_1)$ and $(2C_2)^{-1}\leq R(\tilde{x},0) \leq 2C_2$ for all $\tilde{x}\in \tilde{U}_j$, where $C_1=C_1(n, \theta, \frac{\varepsilon}{2})$ and $C_2(n, \theta, \frac{\varepsilon}{2})$ are constants appearing in Proposition  2.6 in \cite{Hu23} with $\varepsilon$ there replaced by $\frac{\varepsilon}{2}$.  Then for $j$ large enough $\pi_j(\tilde{U}_j)$ is a strong $\varepsilon$-neck or an $\varepsilon$-cap  centered at $(x_j,0)$,  (for the cap case we  need to use  the argument in the proof of Proposition  2.6 in \cite{Hu23},)  with $B_{g^{(j)}(0)}(x_j, (2C_1)^{-1}) \subset U \subset B_{g^{(j)}(0)}(x_j, 2C_1)$, and $(2C_2)^{-1} \leq R(x,0) \leq 2C_2$ for all $x\in U$.  Thus we also get a contradiction.

The derivative estimates follows from Proposition 3.8 (ii).
\hfill{$\Box$}

\section{Existence of $(r, \delta)$-surgical solutions}

The following result extends Proposition 3.1 in \cite{Hu23} to the noncompact case.

\begin{prop} \label{prop 3.1}\ \  (cf. Proposition 3.1 in \cite{Hu23})\ \  Let $0<\varepsilon_0\leq \varepsilon_3$, where the constant $\varepsilon_3$ is as in Proposition 2.4 in \cite{Hu22}.  Fix $0 < \varepsilon <  \tilde{\varepsilon}_1(\varepsilon_0)$, where $\tilde{\varepsilon}_1(\cdot)$ is as in Lemma 2.3 in \cite{Hu22}.
Let $(\mathcal{O},g)$ be a  connected, complete $n$-orbifold with at most isolated singularities and with bounded curvature. Suppose that  each point of $\mathcal{O}$ is a center of an
$\varepsilon$-neck or an $\varepsilon$-cap.  Then $\mathcal{O}$ is diffeomorphic  to a spherical orbifold, or a neck, or a cap, or a connected sum of at most two spherical orbifolds with at most isolated singularities. If we assume in addition  that $\mathcal{O}$ is a manifold, then either $\mathcal{O}$ is diffeomorhic to a spherical space form, or $\mathcal{O}$ is  diffeomorphic  to   a  quotient manifold of $\mathbb{S}^{n-1} \times \mathbb{R}$ by standard isometries.
\end{prop}

\noindent {\bf Proof}.  It remains to consider the case that $(\mathcal{O},g)$ is  noncompact.  In this case, if every point in $\mathcal{O}$ is a center of an $\varepsilon$-neck, then by using Lemma 2.3  and Proposition 2.4 in \cite{Hu22} we get Hamilton's canonical parametrization $\phi: \mathbb{S}^{n-1}/\Gamma \times \mathbb{R}\rightarrow \mathcal{O}$, which is a diffeomorphism, where $\Gamma$ is a finite subgroup of $O(n)$ acting freely on $\mathbb{S}^{n-1}$.    If $\mathcal{O}$ contains an $\varepsilon$-cap, then by using Hamilton's  canonical parametrization again we see that $\mathcal{O}$ is diffeomorphic to $D^n//\Gamma \cup_h \mathbb{S}^{n-1}/\Gamma \times [0, \infty)$, or
 $(\mathbb{S}^n// \langle\Gamma, \hat{\sigma} \rangle \setminus B)\cup_h \mathbb{S}^{n-1}/\Gamma \times [0, \infty)$, or $(\mathbb{S}^n// (x,\pm x') \setminus B) \cup_h  \mathbb{S}^{n-1} \times [0, \infty)$, where $h$ is a diffeomorphism between the boundaries of the relevant manifolds/orbifolds.  By (the proof of) Theorem 2.2 in Chapter 8  of \cite{Hi}, or Proposition 7.7 in \cite{B}, or Proposition 7.6.6 in \cite{Mu}, in this case   $\mathcal{O}$ is diffeomorphic to $\mathbb{R}^n//\Gamma$, or $\mathbb{S}^n// \langle\Gamma, \hat{\sigma} \rangle \setminus \bar{B}$, or $\mathbb{S}^n// (x,\pm x') \setminus \bar{B}$.
     \hfill{$\Box$}

\vspace *{0.2cm}

\noindent {\bf Definition  }(cf. \cite{BBM} and \cite{Hu22}). \ \  A piecewise $C^1$-smooth
  evolving complete Riemannian $n$-orbifold $\{(\mathcal{O}(t), g(t))\}_{t \in I }$ with at most isolated singularities is a
  surgical solution to the Ricci flow if it has the following
  properties.

  i. The equation $\frac{\partial g}{\partial t}=-2 \hspace*{0.4mm} \text{Ric}$ is satisfied
  at all regular times;

  ii. For each singular time $t_0$ there is a   collection
  $\mathcal{S}$ of disjoint embedded $\mathbb{S}^{n-1}/\Gamma$'s in $\mathcal{O}(t_0)$
  (where $\Gamma$'s are finite subgroups of $O(n)$ acting freely on $\mathbb{S}^{n-1})$, and an Riemannian orbifold $\mathcal{O}'$ such that

  (a) $\mathcal{O}'$ is obtained from  $\mathcal{O}(t_0)$ by removing a suitable  open tubular neighborhood of each  element of $\mathcal{S}$ and
  gluing in a Riemannian orbifold  diffeomorphic to $D^n//\Gamma$ along each boundary component thus produced which is diffeomorphic to $\mathbb{S}^{n-1}/\Gamma$;

 (b) $\mathcal{O}_+(t_0)$ is a union of some connected components of $\mathcal{O}'$ and
 $g_+(t_0)=g(t_0)$ on $\mathcal{O}_+(t_0)\cap \mathcal{O}(t_0)$;

(c) each component of $\mathcal{O}'\setminus \mathcal{O}_+(t_0)$ is
diffeomorphic to a spherical orbifold with at most isolated singularities, or a neck, or a cap, or a  connected sum of at most two spherical orbifolds with at most isolated singularities.

\vspace *{0.2cm}

\noindent {\bf Definition} (cf. \cite{P2} and \cite{Hu15}).  Let $\varepsilon$  and $C$ be  positive constants.  we say that a space-time point $(x,t_0)$ in a surgical solution to the Ricci flow has an $(\varepsilon, C)$-canonical neighborhood if  there is   an open
neighborhood $U$ of $x$ satisfying $B(x,t_0,C^{-1}R(x,t_0)^{-\frac{1}{2}}) \subset U\subset
B(x,t_0,CR(x,t_0)^{-\frac{1}{2}})$ and with   $C^{-1}R(x,t_0) \leq R(x',t_0) \leq CR(x,t_0)$ for any $x'\in U$, which falls into one of the following three
types:

(a) $U$ is a  strong $\varepsilon$-neck  centered at $(x,t_0)$,

(b)  $U$ is an $\varepsilon$-cap centered at $(x,t_0)$,

 (c)  $U$ is   compact and strictly PIC2,

\noindent moreover,   the scalar curvature in $U$ at time $t_0$  satisfies
the derivative estimates
\begin{equation*}
|\nabla R|\leq C R^{\frac{3}{2}} \hspace*{8mm} \mbox{and} \hspace*{8mm}
|\frac{\partial R}{\partial t}|\leq C R^2,
\end{equation*}
and in case (c) the sectional curvature in $U$ at time $t_0$ is greater than $C^{-1}R(x,t_0)$.
\vspace *{0.2cm}

\vspace *{0.2cm}

\noindent {\bf Canonical neighborhood assumption}:  Fix positive constants $\varepsilon$ and $C$. Let $r>0$. A surgical solution to the Ricci flow  $\{(\mathcal{O}(t), g(t))\}_{t \in I}$
 satisfies the canonical neighborhood assumption  $(CN)_r$ with $(4\varepsilon, 4C)$-control if  any  space-time point $(x,t)$ with  $R(x,t)\geq
r^{-2}$ has a  $(4\varepsilon, 4C)$-canonical neighborhood.

\vspace *{0.2cm}

We adapt two more definitions from \cite{BBM}.

\vspace *{0.2cm}

\noindent {\bf Definition} (cf. \cite{BBM}, \cite{Hu15} and \cite{Hu22}). \ \ Given an interval
$I\subset [0,+\infty)$, fix surgery parameters $r$, $\delta>0$
and let
$h$, $D$, $\Theta=2Dh^{-2}$ be the associated cutoff parameters  as in \cite{Hu22}. Let
$(\mathcal{O}(t),g(t))$ ($t \in I$) be an evolving complete Riemannian orbifold with  at most isolated singularities. Let $t_0 \in I$ and $(\mathcal{O}_+,g_+)$ be a
(possibly empty) Riemmanian $n$-orbifold. We say that
$(\mathcal{O}_+,g_+)$ is obtained from
$(\mathcal{O}(\cdot),g(\cdot))$ by $(r,\delta)$-surgery at time
$t_0$ if

i. $R_{\text{max}}(g(t_0))=\Theta$, and there is a  collection    of pairwise disjoint strong $\delta$-necks in $\mathcal{O}(t_0)$, centered at some points with scalar curvature equal to $h^{-2}$ at time $t_0$, such that   $\mathcal{O}_+$ is obtained from $\mathcal{O}(t_0)$ by doing
  surgery along these necks, and removing each of the following components:

  (a)  a compact component which is strictly PIC2 and  has sectional curvature bounded below by  $(4C)^{-1}R(x,t_0)$,  where $x$ is some point in this component,

  (b) a component  where each point  is a center of a
$4\varepsilon$-neck or a $4\varepsilon$-cap;

ii. $R_{\text{max}}(g_+)\leq
\Theta/2$ when $\mathcal{O}_+\neq \emptyset$.

\vspace *{0.2cm}

\noindent {\bf Definition} (cf. \cite{BBM} and \cite{Hu22}). \ \  A complete surgical solution
$(\mathcal{O}(\cdot),g(\cdot))$ to the Ricci flow defined on some time interval
$I\subset [0,+\infty)$  is an $(r,\delta)$-surgical solution  if it has the  following
properties:

i.  It satisfies the $(f,\theta)$-pinching assumption  at any time $t\in I$, and $R(x,t) \leq \Theta$ for all $(x,t)$;

ii. At each singular time $t_0\in I$,
$(\mathcal{O}_+(t_0),g_+(t_0))$ is obtained from
$(\mathcal{O}(\cdot),g(\cdot))$ by $(r,\delta)$-surgery at time
$t_0$;

iii. The canonical neighborhood assumption  $(CN)_r$ with $(4\varepsilon, 4C)$-control holds.

\vspace *{0.2cm}

\vspace *{0.2cm}

The following result improves Proposition 4.9 in \cite{Hu22} slightly.

\begin{prop} \label{prop 3.2}(cf. Proposition C in \cite{BBM}, Proposition 10.9 in \cite{B19}, Lemma 4.5 in \cite{CTZ},  Lemma 5.4 in \cite{Hu15}, Proposition 4.9 in \cite{Hu22}, and  Lemma 5.2 in \cite{P2}) Let $(\mathcal{O},g_0)$ be a complete Riemannian orbifold of dimension $n\geq 5$ with at most isolated singularities,  with   isotropic curvature bounded below by a positive constant, with scalar curvature  bounded above, and with $\text{Vol} \hspace*{1mm} B_1(x)\geq v>0$ for every $x \in |\mathcal{O}|$.   Let $\varepsilon$  and $C$ be chosen as in \cite{Hu23}.
Then there exist a positive constant $\kappa$ and a positive function $\tilde{\delta}(\cdot)$ with the following property: If we have an $(r,\delta)$-surgical solution $(\mathcal{O}(t), g(t))$, $t\in [0,T]$, to the Ricci flow with $\delta\leq \tilde{\delta}(r)$ starting with $(\mathcal{O},g_0)$, then  the flow is $\kappa$-noncollapsed  on all scales less than $\varepsilon$.
\end{prop}

\noindent {\bf Proof}.\ \  Consider a space-time point $(x_0, t_0)$ with $|\text{Rm}(\cdot,\cdot)| \leq r_0^{-2}$ on $P(x_0, t_0,r_0, -r_0^2)$, where $r_0 < \frac{r}{C(n, \theta, \varepsilon)}$  for some positive constant $C(n, \theta, \varepsilon)$ to be chosen    (the  case $r_0 \geq \frac{r}{C(n, \theta, \varepsilon)}$ is treated in \cite{Hu22}).
We want to bound $\text{vol}_{t_0}(B(x_0,t_0,r_0))/r_0^n$ from below.  By an argument using orbifold Bishop-Gromov theorem (cf. \cite{B93} and \cite{L}) we may assume w.l.o.g. that there is some point $(x', t') \in \overline{P(x_0, t_0,r_0, -r_0^2)}$ such that $|\text{Rm}(x', t')|=r_0^{-2}$. (Compare Lemma 10.1.2 in  \cite{BBB+}, and p.232-233 in \cite{CZ}.)   We may assume that $R(x',t')\geq r^{-2}$ and the point $x'$ is contained in a  compact component, say $Y$,  of strictly PIC2 at time $t'$  and with sectional curvature at time $t'$ greater than $(4C)^{-1}R(x',t')$, since the other cases are dealt with in \cite{Hu22}.  $Y$ is diffeomorphic to a spherical orbifold $\mathbb{S}^n//\Pi$,  where $\Pi$ is a finite subgroup of $O(n+1)$. So, $Y$ is diffeomorphic to either a spherical manifold $\mathbb{S}^n/\Upsilon$, where $\Upsilon$ is a finite subgroup of $O(n+1)$ acting freely on $\mathbb{S}^n$, or a spherical orbifold $\mathbb{S}^n//\Gamma$,  or   $\mathbb{S}^n// \langle\Gamma, \hat{\sigma} \rangle$, where $\Gamma$ is a finite subgroup of $O(n)$ acting freely on $\mathbb{S}^{n-1}$, cf. \cite{CTZ} and \cite{Hu22}.  To bound $\text{vol}_{t_0}(B(x_0,t_0,r_0))/r_0^n$ from below we only need to bound the topology of $Y$.  So we may assume that  $|\Upsilon|>2$ and
$|\Gamma|>2$.  It suffices to bound the order $|\Pi|$ (or $|\Upsilon|$ and $|\Gamma|$).

Let $t_1=\inf \{t \hspace*{1mm} | \hspace*{1mm}  0\leq t \leq t',  Y\times (t, t']  \hspace*{1mm} \text{is}  \hspace*{1mm} \text{unscathed}\}$.
  Then
  $\inf_{x \in Y} R_{g_+(t_1)}(x)< r^{-2}$, which follows from the definition of  $(r,\delta)$-surgical solution.  Let $t_0'= \sup \{t \hspace*{1mm}  | \hspace*{1mm}  t_1 < t \leq t', \hspace*{1mm}  \inf_{x \in Y} R(x,t)=r^{-2}\}$. Then either $(Y, g(t_0'))$  is a compact component  of strictly PIC2, with $|\text{Rm}| \leq
C(n, \theta, \varepsilon)^2r^{-2}$, and with sectional curvature  greater than $(4C)^{-1}R(\tilde{x},t_0')$ for some $\tilde{x}  \in Y$  with $R(\tilde{x},t_0')=r^{-2}$, or   we  can find a point $x_0'\in Y$  such that $(x_0',t_0')$ is a center of a strong $4\varepsilon$-neck (diffeomorphic to $\mathbb{S}^{n-1}/\Gamma \times \mathbb{R}$) on which  $|\text{Rm}| \leq
C(n, \theta, \varepsilon)^2r^{-2}$. We only need to consider the former case, since the latter case is considered in \cite{Hu22}. Then note that  $B(\tilde{x},t_0',r) \subset Y$, and we have
\begin{equation}
\frac{\text{vol}_{t_0'}(Y)}{r^n}  \geq  \frac{\text{vol}_{t_0'}(B(\tilde{x},t_0',r))}{r^n} \geq c(n, \theta, \varepsilon, g_0)>0
\end{equation}
by (the proof of) Case 1 in the proof of Proposition 4.9 in \cite{Hu22}.  On the other hand,  the Riemannian universal cover $(\tilde{Y}, \tilde{g}(t_0'))$ of  $(Y, g(t_0'))$ has diameter $\leq 2\pi \sqrt{C}r$  by Myers theorem, hence by Bishop volume comparison theorem we have
\begin{equation}
|\Pi| \text{vol}_{t_0'}(Y) = \text{vol}_{t_0'}(\tilde{Y}) \leq   C(n) (2\pi \sqrt{C}r)^n.
\end{equation}
From (3.1) and (3.2) we get a desired upper bound for $|\Pi|$.
Compare  the arguments for Case 2 in Section 10.4 of \cite{BBM} and Step 2 in the proof of Lemma 4.5 in \cite{CTZ}.
\hfill{$\Box$}

\vspace *{0.2cm}

\noindent {\bf Remark}.  In the proof of Lemma 5.4 in \cite{Hu15} we overlooked the possibility for the canonical neighborhood to be a compact component.

\vspace *{0.2cm}

Now we  have

\begin{thm} \label{thm 3.3}(cf. Theorem 5.5 in \cite{Hu15}, and Theorem 4.11 in \cite{Hu22}) \ \  Let $(\mathcal{O},g_0)$ be a complete noncompact
orbifold of dimension $n\geq 12$ with at most isolated singularities, with   isotropic curvature bounded below by a positive constant, with scalar curvature  bounded above, and with $\text{Vol} \hspace*{1mm} B_1(x)\geq v>0$ for every $x \in |\mathcal{O}|$.  Let  $\varepsilon$  and $C$ be chosen as in   \cite{Hu22} and \cite{Hu23}. Then  we can find  positive numbers $\hat{r}$ and $\hat{\delta}$   such that there exists an $(\hat{r}, \hat{\delta})$-surgical solution (with $(4\varepsilon, 4C)$-control) to the Ricci flow starting from $(\mathcal{O},g_0)$, which  becomes extinct in finite time.
\end{thm}

\noindent {\bf Proof}.\ \  With the help of Brendle's curvature pinching estimates, the proof is almost identical to that of  Theorem 4.11 in \cite{Hu22}, except that we need to use Propositions \ref{prop 2.1},  \ref{prop 3.1} and  \ref{prop 3.2} here to replace Propositions 2.7 and 3.1 in \cite{Hu23} and  Proposition 4.9 in \cite{Hu22} repectively,  and something like  Lemma 5.9 in \cite{BBM} to control the difference of two consecutive singular times; cf. also \cite{Hu13}.
 \hfill{$\Box$}

\section {Proof of Theorems 1.1 and 1.2}

\noindent {\bf Proof of Theorem 1.2}.  Let  $(\mathcal{O},g_0)$ satisfy the assumptions of Theorem \ref{thm 1.2}. By Theorem \ref{thm 3.3} we can construct
an $(r,\delta)$-surgical solution to the Ricci flow starting with $(\mathcal{O}, g_0)$,  which becomes extinct in finite time. Recalling that each point in any component that is removed in the process of surgery is contained in a canonical neighborhood, from the proof of Proposition \ref{prop 3.2} we see that the topology of these components  have only finite kinds of diffeomorphism types.
Combining this with  Lemma  5.1 in \cite{Hu22}, which can be extended to  our present situation, we see that the conclusion of Theorem \ref{thm 1.2} holds.
\hfill{$\Box$}

\vspace *{0.2cm}

\noindent {\bf Proof of Theorem 1.1}. Theorem \ref{thm 1.1} follows from Theorem \ref{thm 1.2} and the arguments in the proof of Theorem 1.1 in \cite{Hu23}.
  \hfill{$\Box$}

\vspace *{0.2cm}
From the proof of Theorem \ref{thm 1.1} we have
\begin{thm} \label{thm 4.1}(cf. Theorem 3.3 in \cite{Hu23})   \ \
  Let $(\mathcal{O},g_0)$ be a complete noncompact
orbifold of dimension $n\geq 12$ with at most isolated singularities, with   isotropic curvature bounded below by a positive constant, with scalar curvature  bounded above, and with $\text{Vol} \hspace*{1mm} B_1(x)\geq v>0$ for every $x \in |\mathcal{O}|$. Then   there is a finite collection
$\mathcal{X}$ of spherical $n$-orbifolds  with  at most isolated singularities and   manifolds admitting $\mathbb{S}^{n-1}\times \mathbb{R}$-geometry such that  $\mathcal{O}$ is diffeomorphic to  a
connected sum of  members of $\mathcal{X}$, where the connected sum occurs at regular points.
\end{thm}

\vspace *{0.2cm}

\noindent {\bf Acknowledgements}.  I was partially supported by NSFC (12271040) and Beijing Natural Science Foundation (Z190003).

\vspace *{0.4cm}


\vspace *{0.4cm}

Laboratory of Mathematics and Complex Systems (Ministry of Education),

School of Mathematical Sciences, Beijing Normal University,

Beijing 100875,  People's Republic of China

 E-mail address: hhuang@bnu.edu.cn

\end{document}